\documentclass[12pt,a4paper,reqno]{amsart}
\usepackage{amssymb}
\usepackage{layout} 
\def\boxit{$\sqcap\kern-8pt\sqcup$}
\def\littbox{\null\hfill\boxit{}}

\newtheorem{theorem}{Theorem}
\newtheorem{lemma}{Lemma}

\newcommand{\R}{\mathbb R}

\usepackage{graphics}
\usepackage{pgf,xcolor}
\usepackage{color}
\input epsf

\begin{document}
\title[Edge separators for quasi-binary trees]{Edge separators for quasi-binary trees}
\author[Jorge Luis Ram\'{i}rez Alfons\'{i}n]{Jorge Luis Ram\'{i}rez Alfons\'{i}n}
\address{Institut de Math\'ematiques et de Mod\'elisation de Montpellier, Universit\'e Montpellier 2, 
Place Eug\`ene Bataillon, 34095 Montpellier}
\email{jramirez@math.univ-montp2.fr}
\author[Serge  Tishchenko
]{Serge  Tishchenko}
\address{IMJ, Universit\'{e} Pierre et Marie Curie,
Paris 6, 4 Place Jussieu, 75252 Paris Cedex 05} 

\begin{abstract}  One wishes to remove $k-1$ edges of a vertex-weighted tree $T$ such
that the weights of the $k$ induced connected components are approximately the same.  How well can one do it ?
In this paper, we investigate such $k$-separator for {\em quasi-binary} trees.  We show that, 
under certain conditions on the total weight of the tree, a particular $k$-separator can be constructed such that the smallest (respectively the largest) weighted component is lower (respectively upper) bounded.  Examples showing optimality for the lower bound are also given.
\end{abstract}

\maketitle
\noindent
\textbf{Keywords:} Binary tree, separator
\smallskip


\section{Introduction} \label{sec:intro}

The seminal paper by Lipton and Tarjan \cite{LT} has inspired a number of separator-type problems and applications (we refer the reader to \cite{RH} for a recent survey on separators). 

Let us consider the following question. 

\begin{quote} One wishes to split a given embedding of a planar connected graph $G$ into 
blocks formed by weighted faces (weights might be thought as area of faces) such that  
the dual of the planar graph induced by each block is connected and the blocks' weights are approximately the same.    
How well can this be done ?
\end{quote} 

One way to answer the latter  is by considering $k$-separators on a spanning tree 
$T_G$ of the vertex-weighted dual graph of $G$.
Indeed, one may want to remove $k-1$ edges of $T_G$ such
that the weights of the $k$ induced connected components of $T_G$ are approximately the same.  

More formally, let $T=(V,E)$ be a graph and let $\omega : V(T)\longrightarrow \R$ be a {\em weight} fonction.
Let $\omega(T)=\sum\limits_{v\in V(T)}\omega(v)$ and
let $1\le k\le |V|-1$ be an integer. A  {\em $k$-separator} of $T$ is a set $F\subset E(T)$ with $|F|=k-1$ 
whose deletion induce $k$ connected components, say $C_{F}^1(T),\dots, C_{F}^{k}(T)$. 
If we let $\omega(C_{F}^{i}(T))=\sum\limits_{v\in V(C_{F}^{i}(T))}\omega(v)$ then  
$\omega(T)=\sum\limits_{i=1}^k\omega(C_{F}^{i}(T)$.   Let
 
$$\beta_k(T):=\max\limits_{F\subseteq E, |F|=k-1}\left\{\min\limits_{1\le i\le k}\omega(C_{F}^{i}(T))\right\}$$ and 

$$\alpha_k(T):=\min\limits_{F\subseteq E, |F|=k-1}\left\{\max\limits_{1\le i\le k}\omega(C_{F}^{i}(T))\right\}.$$

An optimal $k$-separator is achieved when $\beta_k(T)=\alpha_k(T)={1\over k}\omega(T)$.

In this paper, we investigate the existence of $k$-separators with large (resp. small) value for $\beta_k$ (resp. for $\alpha_k$) for the class of {\em quasi-binary} trees.
A tree is called {\em binary} if the degree of any vertex is equals three except for {\em pending} vertices (vertices of degree one) and a {\em root} vertex (vertex of degree two).
A tree is say to be {\em quasi-binary} if it is  a connected subgraph of a binary tree. Notice that a good $k$-separators in quasi-binary trees would lead to a good $k$-block separators for triangulated planar graphs in the above question.

Since $d(v)=1,2$ or $3$ for any $v\in V(T)$ of a quasi-binary tree $T$ then we may define, for each  $i=1,2,3$,

$$V_i:=\{v\in V(T) | d(v)=i\}\ \ \text{ and } \ \ \omega_i:=\max\{\omega(v) | v\in V_j\ \text{for each }\  i\le j\le 3\}.$$ 

Notice that

\begin{equation}\label{eq1}
V(T)=V_1\cup V_2\cup V_3,\ \  \omega_1\ge \omega_2\ge \omega_3\ \ \text { and }\ \ \omega_1n_1+\omega_2n_2+\omega_3n_3\ge \omega(T). 
\end{equation}
where $n_i=|V_i|$ for each $i=1,2,3$.

Our main results are the following.

\begin{theorem}\label{mt1} Let $T$ be a quasi-binary tree. Let $k\ge 2$ be an integer and $\gamma\in\R$ with $\gamma\ge \omega_3$. If 
$$\omega(T)\ge \max\left\{{{(k+1)(k-2)}\over {(k-1)}}\omega_1-{2\over {(k-1)}}\gamma,{{2(k+1)(k-2)}\over {(k-1)}}\omega_2-k\gamma\right\}$$ then 
$$\alpha_k(T)\le {{2\omega(T)+(k-1)\gamma}\over {k+1}}\cdot$$
\end{theorem}

\begin{theorem}\label{mt2} Let $T$ be a quasi-binary tree. Let $k\ge 2$ be an integer and $\gamma\in\R$ with $\gamma\ge \omega_3$. If 
$$\omega(T)\ge \max\left\{{{(2k-1)}\over {2}}\omega_1-{1\over 2}\gamma,(2k-1)\omega_2-k\gamma\right\}$$ then 
$$\beta_k(T)\ge {{\omega(T)-(k-1)\gamma}\over {2k-1}}\cdot$$
\end{theorem}

We notice that the bounds for $\alpha_k(T)$ and $\beta_k(T)$ are not necessarily reached by using the same $k$-separator.

The second author has studied $k$-separators in a more general setting (for planar graphs with weights on vertices, edges and faces) where a lower bound for $\beta_k$  is obtained \cite{T1} . 
We noticed that the conditions given in \cite{T1} are different from  those presented in Theorem \ref{mt2}. 
Moreover, the proof of Theorem \ref{mt2} (which is in the same spirit as that of Theorem \ref{mt1})
is different from that given in \cite{T1}. The value $\alpha_k$ is not treated in \cite{T1} at all.  
\smallskip

In the following section we present some preliminary results needed for the rest of the paper. 
Main results are proved in Section \ref{mainsec}. Finally, a family of  quasi-binary trees, 
showing optimality of Theorem \ref{mt2}, is constructed in the last section.  

\section{Preliminary results}

Let $T$ be a quasi-binary tree with $\omega(T)>0$. We let $n=|V(T)|>1$, and $n_i=|V_i|$ for each $i=1,2$ and $3$.
We observe that

\begin{equation}\label{eq2}
n_1+2n_2+3n_3=2|E(T)|=2(n-1)=2n_1+2n_2+2n_3-2 \ \text{and thus}\  n_1=n_3+2.
\end{equation}
 
Our main theorems will be proved by induction. For, we need the following two lemmas.

\begin{lemma}\label{l1} Let $T$ be a quasi-binary tree with $n=|V(T)|> 1$. Let $\gamma,\eta \in\R$ such that $\gamma\ge \omega_3$ and $\max\left\{{{\omega_1-\gamma}\over {2}},\omega_2-\gamma\right\}\le \eta \le {{\omega(T)}\over 2}$.
Then, there exist $e\in E(T)$ such that 
$$\eta \le \omega(C_e^i)\le 2\eta+\gamma$$ 
for some $1\le i\le 2$, where $C_e^1,C_e^2$ denote the two connected components of $T\setminus \{e\}$ .
\end{lemma}

{\em Proof.} The inequality $\eta \le \omega(C_e^i)$ holds for $i=1$ and/or $i=2$ and for any $e\in E(T)$, otherwise $\omega(T)={{\omega(T)}\over 2}+{{\omega(T)}\over 2}\ge 2\eta>\omega(C_e^1)+\omega(C_e^2)=\omega(T)$, which is a contradiction. 

We now prove the right-hand side inequality. Without lost of generality, we suppose that  
$\omega(C_e^1)\ge \eta$ for each $e\in E(T)$. If we also have that $\omega(C_e^1)\ge \eta$ then
we choose indices such that $|V(C_e^1)|\le |V(C_e^2)|$.  
 
We proceed by contradiction, suppose that $\omega(C_e^1)>2\eta+\gamma$ for all $e\in E(T)$.  Let $e=\{v_1,v_2\}$ with $v_i\in V(C_e^i)$ be the edge that minimizes $|V(C_e^1)|$. We have three cases. 
\smallskip

Case 1) If $d(v_1)=1$ then $\omega(C_e^1)=\omega(v_1)\le \omega_1$. Since $\eta\ge {{\omega_1-\gamma}\over 2}$ then $2\eta+\gamma\ge\omega_1\ge\omega(C_e^1)$, which is a contradiction.
\smallskip

Case 2) If $d(v_1)=2$ then we let $f=\{v_1,u\}\in E(T), f\neq e$ be the other edge incident to $v_1$. 

Let $C_f^i$, $i=1,2$ be the two connected components of $T\setminus \{f\}$. 
Since $|V(C_f^1)|\ge |V(C_e^1)|$ then $V(C_e^1))=V(C_f^2)\cup\{v_1\}$ so 
$\omega(C_f^2)=\omega(C_e^1)-\omega(v_1)>  2\eta+\gamma-\omega_2\ge \eta$, 
and thus $|V(C_f^2)|\ge |V(C_f^1)|\ge |V(C_e^1)|$  which is a contradiction. 
\smallskip

Case 3) If $d(v_1)=3$ then we let $f_1=\{v_1,u\},f_2=\{v_1,v\}\in E(T), f_1,f_2\neq e$ be the other two edges incident to $v_1$
with $V(C_{f_1}^2)\cup V(C_{f_2}^2)=V(C_e^1)\setminus \{v_1\}$. So,
$\omega(C_{f_1}^2)+\omega(C_{f_2}^2)=\omega(C_e^1)-\omega(v_1)>  2\eta+\gamma-\omega_3\ge 2\eta$. Without loss of generality, we suppose that $\omega(C_{f_1}^2)\ge \omega(C_{f_2}^2)$, and thus
$\omega(C_{f_1}^2)>\eta$ and $|V(C_{f_1}^2)|\ge |V(C_{f_1}^1)|\ge |V(C_{e}^1)|$ which is a contradiction.
\littbox

\begin{lemma}\label{l2} Let $T$ be a quasi-binary tree with $n=|V(T)|> 1$. Let $\gamma \in\R$ such that $\gamma\ge \omega_3$. If $\omega(T)\ge \max\left\{{{3\omega_1-\gamma}\over {2}},3\omega_2-2\gamma\right\}$ then

$$\beta_2(T)\ge {{\omega(T)-\gamma}\over {3}} \ \ \text {and} \ \   \alpha_2(T)\le {{2\omega(T)+\gamma}\over {3}}\cdot$$
\end{lemma}

{\em Proof.} We first claim that $\omega(T)\ge -2\gamma$. Indeed,

$$\begin{array}{ll}
(n+2)(\omega(T) +2\gamma)& = (2n_1+n_2)\omega(T)+(2n_1+n_2)2\gamma\\
& \ge 2n_1\left({3\over 2}\omega_1-{1\over 2}\gamma\right)+n_2(3\omega_2-2\gamma)+3n_3(\omega_3-\gamma)+(2n_1+n_2)2\gamma\\
& \ge 3(n_1\omega_1+n_2\omega_2+n_3\omega_3)+3(n_1-n_3)\gamma\ge 3(\omega(T)+2\gamma).\\ 
\end{array}$$

Therefore, since $n\ge 1$ then $(\omega(T)+2\gamma)\ge 0$ and the result follows. 

We now claim that  $$\max\left\{{{\omega_1-\gamma}\over {2}},\omega_2-\gamma\right\}\le \eta \le {{\omega(T)}\over 2}$$ 
is verified by taking $\eta={{\omega(T)-\gamma}\over 3}$. Indeed, 
$$\eta\le  {{\omega(T)}\over 2} \iff {{\omega(T)-\gamma}\over 3}\le  {{\omega(T)}\over 2} \iff \omega(T)\ge -2\gamma.$$

For the lower bound, we have two cases.

Case 1)  
$$\eta\ge{{\omega_1-\gamma}\over 2} \iff {{\omega(T)-\gamma}\over 3} \ge {{\omega_1-\gamma}\over 2} \iff \omega(T)\ge {{3\omega_1-\gamma}\over 2}$$
which is true by hypothesis.

Case 2)
$$\eta\ge \omega_2-\gamma \iff {{\omega(T)-\gamma}\over 3} \ge \omega_2-\gamma \iff \omega(T)\ge 3\omega_2-2\gamma$$
which is true by hypothesis.

Therefore, by Lemma \ref{l1}, there is an edge $e\in E(T)$ such that   
 $${{\omega(T)-\gamma}\over {3}} \le \omega(C_e^i)  \le {{2\omega(T)+\gamma}\over {3}}$$
for one of the two connected components $C_e^i, \ i=1,2$ of $T\setminus \{e\}$
and the result follows.
\littbox

\section{Proofs of main results}\label{mainsec}

We may now prove our main results.

{\em Proof of Theorem \ref{mt1}.} We first show that $\omega(T)>-k\gamma$ (needed for the rest of the proof). 
For, we consider

$$\begin{array}{ll} 
\left(n-\frac{2(k-1)}{k}\right)(\omega(T)+k\gamma) & =(2n_1+n_2)\omega(T)+(2n_1+n_2)k\gamma-{{2(2k-1)}\over k}(\omega(T)+k\gamma)\\
& \ge 2n_1\left({{(2k-1)}\over 2}\omega_1-{1\over 2}\gamma\right)+n_2((2k-1)\omega_2-k\gamma)\\
&\ \ \ +(2k-1)n_3(\omega_3-\gamma)+(2n_1+n_2)k\gamma-{{2(2k-1)}\over k}(\omega(T)+k\gamma)\\
& \ge (2k-1)(n_1\omega_1+n_2\omega_2+n_3\omega_3)+(2k-1)(n_1-n_3)\gamma\\
&\ \ \ -{{2(2k-1)}\over k}\omega(T)-2(2k-1)k\gamma\\
&\ge {{(2k-1)(k-2)}\over k}\omega(T)>0.\\
\end{array}$$

Since $n>1$ if and only if $n-{{(2k-1)}\over k}>0$ then $(\omega(T)+k\gamma)>0$ and the inequality follows. 

We now shall construct the desired $k$-separator as follows. Let $T_k=T$, we first find an edge $e_k\in E(T_k)$ 
(by using Lemma \ref{l1}) such that one of the connected components of $T_k\setminus \{e_k\}$, say $T_{k-1}$, has a {\em suitable} weight  (the other connected component of  $T_k\setminus \{e_k\}$, say $R_{k-1}$, remains  fixed for the rest of the construction).
By a suitable weight we mean a weight such that Lemma \ref{l1} can be applied to $T_{k-1}$ in order to find 
an edge $e_{k-1}\in E(T_{k-1})$ such that one of the connected components of $T_{k-1}\setminus \{e_{k-1}\}$, say $T_{k-2}$, has again a  suitable weight  (and again the other connected component of  $T_{k-1}\setminus \{e_{k-1}\}$, say $R_{k-2}$, remains  fixed for the rest of the construction), and so on. We claim that the weight of component $T_j$ is suitable if

\begin{eqnarray}\label{bb21}
\frac{(j-1)(k-1)}{(k+1)(k-2)}\omega(T)+\left(\frac{2(j-1)}{(k+1)(k-2)}-1\right)\gamma\le \omega(T_j)\le \frac{(j+1)}{(k+1)}\omega(T)+\frac{(k-j)}{(k+1)}\gamma\cdot
\end{eqnarray}

Now, in order to apply Lemma \ref{l1} we need to define an appropriate parameter $\eta_j$  (that ensures suitable weights throughout the construction). For each $j=k,k-1\dots ,2$, we set 

\begin{eqnarray}\label{bb22}
\eta_{j}=\frac{(k-3)}{2(2k-j-3)}\omega(T_j)-\frac{(j-3)(k-1)}{2(2k-j-3)(k+1)}\omega(T)-\frac{(k+3)(k-2)-j(k-1)}{2(2k-j-3)(k+1)}\gamma\cdot
\end{eqnarray}

We first claim that 

$$\max\left\{\frac{\omega_1-\gamma}{2},\omega_2-\gamma\right\} \le \eta_j\le \frac{\omega(T_j)}{2}\cdot$$

For the lower bound we consider the following

$$\begin{array}{ll}
n_j & = \frac{(k-3)}{2(2k-j-3)}\omega(T_j) - \frac{(j-3)(k-1)}{2(2k-j-3)(k+1)}\omega(T)-\frac{(k-3)(k-2)-j(k-1)}{2(2k-j-3)(k+1)}\gamma\\
 & \ge \frac{(k-3)}{2(2k-j-3)}\left(\frac{(j-1)(k-1)}{(k+1)(k-2))}\right)\omega(T)+
\frac{(k-3)}{2(2k-j-3)}\left(\frac{2(j-3)}{(k+1)(k-2)}-1\right)\gamma\\
& -\frac{(j-3)(k-1)}{2(2k-j-3)(k+1)}\omega(T)-\frac{(k-3)(k-2)-j(k-1)}{2(2k-j-3)(k+1)}\gamma\\
&= \frac{(k-1)}{2(k+1)(k-2)}\omega(T)-\frac{k^2-k-4}{2(k+1)(k-2)}\gamma\\
& \ge \max\left\{\frac{\omega_1}{2}-\frac{1}{(k+1)(k-2)}\gamma,\omega_2-\frac{k(k-1)}{2(k+1)(k-2)}\gamma\right\}-\frac{(k^2-k-4)}{2(k+1)(k-2)}\gamma\\
&=\max\left\{\frac{\omega_1}{2}-\frac{\gamma}{2},\omega_2-\gamma\right\}.
\end{array}$$

And, for  the upper bound, we consider the following. 

$$\begin{array}{ll}
n_j & = \frac{(k-3)}{2(2k-j-3)}\omega(T_j) - \frac{(j-3)(k-1)}{2(2k-j-3)(k+1)}\omega(T)-\frac{(k-3)(k-2)-j(k-1)}{2(2k-j-3)(k+1)}\gamma\\
& =  \frac{\omega(T_j)}{2}-\frac{(k-j)}{2(2k-j-3)}\omega(T_j)- \frac{(j-3)(k-1)}{2(2k-j-3)(k+1)}\omega(T)- \frac{(k+3)(k-2)-j(k-1)}{2(2k-j-3)(k+1)} \gamma\\
& \le \frac{\omega(T_j)}{2}-\frac{(k-j)}{2(2k-j-3)}\left( \frac{(j-1)(k-1)}{(k+1)(k-2)}\omega(T)+\left( \frac{2(j-1)}{(k+1)(k-2)}-1\right)\gamma\right)\\
& \ \ \ - \frac{(j-3)(k-1)}{2(k-j-3)(k+1)}\omega(T)-\frac{(k+3)(k-2)-j(k-1)}{2(2k-j-3)(k+1)}\gamma\\
& \le \frac{\omega(T_j)}{2}-\frac{(j-2)}{(k+1)(k-2)}\left( \frac{(k-1)}{2}\omega(T)+\gamma\right)\\
& =  \frac{\omega(T_j)}{2}-\frac{(j-2)}{2k}\omega(T)-\frac{(j-2)}{k(k+1)(k-2)}\left(\omega(T)+k\gamma\right)\le
\frac{\omega(T_j)}{2}\cdot\\
\end{array}$$

Therefore, by Lemma \ref{l1}, one of the connected components of $T_j\setminus \{e_j\}$, say $R_{j-1}$,  verifies 

\begin{eqnarray}\label{bb25}
\eta_j\le \omega(R_{j-1})\le 2\eta_j+\gamma
\end{eqnarray}

and thus,  the weight of the other connected component of $T_j\setminus \{e_j\}$, says $T_{j-1}$, satisfies

$$\omega(T_j)-2\eta_j-\gamma\le \omega(T_{j-1})\le\omega(T_j)-\eta_j\cdot$$

So, the set of edges $e_1,\dots ,e_{k-1}$ chosen as above gives a $k$-separator $T$ where
the connected component with the largest weight is given by
$\max\limits_{1\le i \le k}\left\{ \omega(T_j) \right\}$. In order to upper bound the latter, we shall show that
that the components $T_j$, $j=k,k-1,\dots ,1$ have suitable weights satisfying both inequalities of \eqref{bb21}.  

We proceed by induction on $j$. If $j=k$ the upper bound is immediate. For the lower bound we have,

\begin{eqnarray}\label{bb23}
\omega(T_k) =\omega(T) & =\frac{(k-1)(k-1)}{(k+1)(k-2)}\omega(T)+\left(\frac{2(k-1)}{(k+1)(k-2)}-1\right)\gamma+\frac{(k-3)}{(k+1)(k-2)}(\omega(T)+k\gamma)\nonumber\\
 & \ge \frac{(k-1)(k-1)}{(k+1)(k-2)}\omega(T)+\left(\frac{2(k-1)}{(k+1)(k-2)}-1\right)\gamma.
\end{eqnarray}

The latter inequality uses the fact that $\omega(T)>-k\gamma$. Suppose that inequalities hold for $j\le k$. By using \eqref{bb21},\eqref{bb22} and \eqref{bb23}, we have

$$\begin{array}{ll}
\omega(R_{j-1}) & \ge \omega(R_j)-2\eta_j-\gamma\\
& = \omega(R_j)-\frac{(k-3)}{(2k-j-3)}\omega(R_j)+\frac{(j-3)(k-1)}{(2k-j-3)(k+1)}\omega(T)+\frac{(k-3)(k-2)-j(k-1)}{(2k-j-3)(k+1)}\gamma-\gamma\\
& \ge \frac{(k-j)}{(2k-j-3)}\left(\frac{(j-1)(k-1)}{(k+1)(k-2)}\right)\omega(T)+\frac{(k-j)}{(2k-j-3)}\left(\frac{2(j-1)}{(k+1)(k-2)}-1\right)\gamma\\
& \ \ \ +\frac{(j-3)(k-1)}{(2k-j-3)(k+1)}\omega(T)-\frac{(k+3)(k-2)-j(k-1)}{(2k-j-3)(k+1)}\gamma-\gamma\\
& = \frac{(j-3)(k-1)}{(k+1)(k-2)}\omega(T)+\left(\frac{2(j-2)}{(k+1)(k-2)}-1\right)\gamma.
\end{array}$$

And

$$\begin{array}{ll}
\omega(R_{j-1}) & \le \omega(R_j)-\eta_j\\
& = \omega(R_j)-\frac{(k-3)}{2(2k-j-3)}\omega(R_j)+\frac{(j-3)(k-1)}{2(2k-j-3)(k+1)}\omega(T)+\frac{(k-3)(k-2)-j(k-1)}{(2k-j-3)(k+1)}
\gamma\\
& \le  \frac{(3k-2j-3)(j+1)}{2(2k-j-3)(k+1)}\omega(T)+\frac{(3k-2j-3)(k-j)}{2(2k-j-3)(k+1)}\gamma\\
& \ \ \ +\frac{(j-3)(k-1)}{2(2k-j-3)(k+1)}\omega(T)+\frac{(k+3)(k-2)-j(k-1)}{2(2k-j-3)(k+1)}\gamma\\
& =\frac{j}{(k+1)}\omega(T)+\frac{(k-j+1)}{(k+1)}\gamma.\\
\end{array}$$

Therefore, \eqref{bb21} holds for all $j=k,\dots ,1$ when $T$ is decomposed into the $k$ components $T_1, R_k,\dots ,R_2$.  So,

$$\begin{array}{ll} 
\alpha_k & = \max\left\{\max\limits_{2\le j \le k}\left\{\omega(R_j) \right\},\omega(T_1) \right\}\overset{\eqref{bb25}}{\le}
\max\left\{\max\limits_{2\le j \le k}\left\{2\eta_j+\gamma \right\},\omega(T_1) \right\}\\
& \overset{\eqref{bb22}}{=} \max\left\{\max\limits_{2\le j \le k}\left\{\frac{(k-3)}{(2k-j-3)}\omega(T_j)-\frac{(j-3)(k-1)}{(2k-j-3)(k+1)}\omega(T)+\frac{k^2-2k+3-2j}{(2k-j-3)(k+1)}\gamma\right\},\omega(T_1)\right\}\\
& \overset{\eqref{bb21}}{\le} \max\left\{\max\limits_{2\le j \le k}\left\{\frac{2}{(k+1)}\omega(T)+\frac{(k-1)}{(k+1)}\gamma\right\},\frac{2}{(k+1)}\omega(T)+\frac{(k-1)}{(k+1)}\gamma\right\}\\
&=\frac{2}{(k+1)}\omega(T)+\frac{k}{(k+1)}\gamma,
\end{array}$$

as desired.
\littbox

{\em Proof of Theorem \ref{mt2}.} We first show that $\omega(T)>-k\gamma$ (needed for the rest of the proof). 
For, we consider

$$\begin{array}{ll} 
\left(n-\frac{2(k-1)}{k}\right)(\omega(T)+k\gamma) 
& = (n+2)(\omega(T)+k\gamma)-\frac{2(2k-1)}{k}(\omega(T)+k\gamma)\\

& =(2n_1+n_2)\omega(T)+(2n_1+n_2)k\gamma-{{2(2k-1)}\over k}\omega(T)-2(2k-1)\gamma\\

& \ge 2n_1\left({{(2k-1)}\over 2}\omega_1-{1\over 2}\gamma\right)+n_2((2k-1)\omega_2-k\gamma)\\
&\ \ \ +(2k-1)n_3(\omega_3-\gamma)+(2n_1+n_2)k\gamma-{{2(2k-1)}\over k}\omega(T)-2(2k-1)\gamma\\

& \ge (2k-1)(n_3\omega_3+n_2\omega_2+n_1\omega_1)-{{2(2k-1)}\over k}\gamma(T)\\
& \ \ \ +(2k-1)(n_1-n_3-2)\gamma\\
&\ge (2k-1)\omega(T)-{{2(2k-1)}\over k}\omega(T)=\frac{(2k-1)(k-2)}{k}\omega(T)>0.\\
\end{array}$$

Since $n>1$ if and only if $n-{{(2k-1)}\over k}>0$ then $(\omega(T)+k\gamma)>0$ and the inequality follows. 

We shall construct the desired $k$-separator in a similar way as done in Theorem \ref{mt1}.
Let $T_k=T$, we find an edge $e_k\in E(T_k)$ (by using Lemma \ref{l1}) such that one of the connected components of $T_k\setminus \{e_k\}$, say $R_k$, has a prescribed weight and which will be fixed for the rest of the construction. By applying Lemma \ref{l1} to the other component of $T_k\setminus \{e_k\}$, say $T_{k-1}$, we find an edge $e_{k-1}\in E(T_{k-1})$ such that one of the connected components of $T_{k-1}\setminus \{e_{k-1}\}$, say $R_{k-2}$, has a prescribed weight and which will be fixed for the rest of the construction, and so on. The only difference with the procedure in the proof of Theorem \ref{mt1}
is that the value $\eta_j$ is now fixed for any step of the construction
$$\hbox{$\eta_j=\eta=\frac{\omega(T)-(k-1)\gamma}{2k-1}$ for all $j=k,k-1,\dots ,2$.}$$

First, we claim that $\eta\ge\max\left\{\frac{\omega_1-\gamma}{2},\omega_2-\gamma\right\}$. Indeed, 
$$\begin{array}{ll}
\eta &= \frac{\omega(T)-(k-1)\gamma}{2k-1}\\
& \ge \max\left\{\frac{(2k-1)\omega_1-\gamma}{2(2k-1)},\frac{(2k-1)\omega_2-k\gamma}{2(2k-1)}\right\}-\frac{(k-1)\gamma}{2k-1}\\
& = \max\left\{\frac{\omega_1-\gamma}{2},\omega_2-\gamma\right\}.\\
\end{array}$$

Therefore, at each step (by Lemma \ref{l1}) one of the connected components of $T_j\setminus\{e_j\}$, say $R_{j-1}$ verifies

$$\eta=\frac{\omega(T)-(k-1)\gamma}{2k-1}\le \omega(R_{j-1})\le \frac{2\omega(T)+\gamma}{2k-1}=2\eta+\gamma.$$

The weight of the other connected component of $T_j\setminus\{e_j\}$, say $T_{j-1}$ satisfies
$$\omega(T_j)-\frac{2\omega(T)-\omega}{2k-1}\le\omega(T_{j-1})\le\omega(T_j)-\frac{\omega(T)-(k-1)\gamma}{2k-1}.$$

Since $\omega(T_k)=\omega(T)$, we obtain
$$\omega(T_j)\ge \omega(T)-(k-j)\frac{2\omega(T)+\gamma}{2k-1}=\frac{(2j-1)\omega(T)-(k-j)\gamma}{2k-1}, \ \ j=k,k-1,\dots , 1.$$

We claim $\eta\le \frac{\omega(T_j)}{2}$ for each  $j=k,k-1,\dots , 2$. Indeed,

$$\begin{array}{ll}
\frac{\omega(T_j)}{2} & \ge \frac{(2j-1)\omega(T)-(k-j)\gamma}{2(2k-1)}\\
& = \frac{\omega(T)-(k-1)\gamma}{2k-1}+\frac{(2j-3)\omega(T)+(k+j-2)\gamma}{2(2k-1)}\\
& =\frac{\omega(T)-(k-1)\gamma}{2k-1}+\frac{(j-2)}{2k}\omega(T)+\frac{(k+j-2)}{2k(2k-1)}(\omega(T)+k\gamma)\\
& \ge \frac{\omega(T)-(k-1)\gamma}{2k-1}=\eta.
\end{array} $$

So, the set of edges $\{e_k,e_{k-1},\dots ,e_2\}$ chosen as above forms a $k$-separator $S_k$ of $T$ where the connected component with the smallest weight is given by

$$\beta(S_k)=\min\{\omega(R_{k-1}),\omega(R_{k-2}),\dots ,\omega(R_{1}),\omega(T_1)\}\ge \frac{\omega(T)-(k-1)\omega}{2k-1}$$
as desired.
\littbox

\section{Tightness}

In this section we show that the lower bound presented in Theorem \ref{mt2} is optimal. For, we consider the quasi-binary tree $T_{k}$ consisting of a root vertex $r$ joined by $k-1$ different paths to $k-1$ vertices $x_1,\dots ,x_{k-1}$ each of which is adjacent to exactly  two vertices of degree one.


We set $\omega(x_i)=\omega>0$ for all $i$, $\omega(r)=\omega(v)=\omega'\ge \omega>0$ where $d(v)=1$
and the weight of any other vertex equals zero.  So, 
$$\omega(T_k)=(k-1)\omega+2(k-1)\omega'+\omega'=(k-1)\omega+(2k-1)\omega'.$$
Let $F$ be an optimal $k$-separator of $T$.  We have that either $F$ contains one of the edges $\{x_i,v\}$, $1\le i\le k-1$ 
with $v$ a pending vertex (so vertex $v$ will be a connected component itself in the separator and thus $\beta_k=\omega'$) or $F$ contains no such edges in which case we find (by an easy analysis of $T_k$) that the root vertex $r$ will be in a connected component containing just vertices of weight zero in any optimal separator (obtaining again that $\beta_k=\omega'$).  
\smallskip

Lower bound of Theorem \ref{mt2} gives 
$$\beta_k\ge \frac{1}{2k-1}\omega(T_k)-\left(\frac{k-1}{2k-1}\right)\omega_3=\frac{1}{2k-1}\left((k-1)\omega+(2k-1)\omega'\right)-\left(\frac{k-1}{2k-1}\right)\omega=\omega'$$
showing the desired optimality.

\end{document}